\input amstex.tex
\input amsppt.sty

\define\E{\Bbb E}
\define\N{\Bbb N}
\define\R{\Bbb R}
\redefine\H{\Bbb H}
\define\C{\Bbb C}
\define\Z{\Bbb Z}
\define\tr{\operatorname{Tr}}

\NoBlackBoxes

\topmatter
\title CLT for spectra of submatrices of Wigner random matrices
\endtitle

\author Alexei Borodin
\endauthor

\abstract We prove a CLT for spectra of submatrices of real symmetric and
Hermitian Wigner matrices. 
We show that if in the standard normalization 
the fourth moment of the off-digonal entries is GOE/GUE-like then
the limiting Gaussian process can be viewed as a collection of simply 
yet nontrivially correlated two-dimensional
Gaussian Free Fields.
\endabstract
\endtopmatter

\document

\subhead Introduction
\endsubhead
Gaussian global fluctuations of eigenvalues of GUE, GOE, Wigner random matrices,
and their generalizations is a well-studied subject, see e.g. Chapter 2 of
\cite{AGZ}
and Chapter 9 of \cite{BS} as well as references therein. One would usually
concentrate on studying the spectrum of the full matrix, but it comes as no
surprise that for large submatrices with a regular limiting behavior,
the joint fluctuations would still be Gaussian. We prove this fact by a slight
modification of the moment method presented in \cite{AGZ}.

It becomes more interesting when one looks at the limiting covariance structure.
In what follows we assume that in the standard normalization 
the fourth moment of the off-diagonal entries of our matrices is the same as for GOE/GUE. 

The first statement is that for such a (real symmetric or Hermitian) Wigner
matrix, the joint fluctuations of spectra of nested submatrices formed by
cutting out top left corners are described by the two-dimensional Gaussian Free
Field (GFF), see e.g. \cite{S} for definitions and basic properties of GFFs.

Although this result seems to be new, the appearance of the GFF is also not
too surprising. Indeed, as was shown in \cite{JN} and
\cite{OR}, for GUE the eigenvalue ensemble of nested matrices arises
as a limit of random surfaces, and for random surfaces the relevance of the GFF
is widely anticipated, see \cite{K}, \cite{BF} for rigorous results and
further references. One might argue however that the GFF interpretation 
simplifies the description of the covariance in the one-matrix case, cf. 
Proposition 3 below. 

The real novelty comes when one considers joint fluctuations for different nested
sequences of submatrices. For each of the nested sequences the fluctuations are again
described by the GFF. On the other hand, when different sequences have
nontrivial and asymptotically regular intersections, these GFFs are correlated,
and the exact form of the covariance kernel turns out to be simple. One could
argue that it is as simple as one could hope for.

The resulting Gaussian process unites a large family of mutually
correlated GFFs. Even for two GFFs the resulting Gaussian process seems to 
be new. An efficient description of the largest natural state space for this
Gaussian process remains an open problem.

It is natural to ask how univeral the limiting process is. We believe that
it also
arises in the world of random surfaces, although it is not {\it a priori\/} clear how
to vary the nested sequence there. The answer comes from representation theory ---
one views random surfaces as originating from restricting suitable representations
to a maximal commutative subalgebra and then one varies that subalgebra. We will
address these models in a later publication.
\smallskip
\noindent{\it Acknowledgements.\/} The author is very grateful to Grigori Olshanski and 
Ofer Zeitouni for valuable comments. The work was partially supported by NSF grant
DMS-1056390.

\subhead Wigner matrices
\endsubhead
Let $\{Z_{ij}\}_{j> i\ge 1}$ and $\{Y_i\}_{i\ge 1}$ be two families of independent
identically distributed real-valued random variables with zero mean such that for any
$k\ge 1$
$$
\max(\E |Z_{12}|^k,\E|Y_1|^k)<\infty.
$$

Assume also that
$$
\E Y_1^2=2,\qquad \E Z_{12}^2=1,\qquad \E Z_{12}^4=3.
$$

Define a (real symmetric) {\it Wigner matrix\/} $X$ by
$$
X(i,j)=X(j,i)=\cases Z_{ij},& i<j,\\ Y_i,&i=j.\endcases
$$

An Hermitian variation of the same definiton is as follows:
Let $\{Z_{ij}\}_{j> i\ge 1}$ now be complex-valued (i.i.d. mean zero) random
variables with the same uniform bound on all
moments. Assume that
$$
\E Y_1^2=1, \qquad \E |Z_{12}|^2=1,\qquad \E |Z_{12}|^4=2.
$$

Define an {\it Hermitian\/} Wigner matrix $X$ by
$$
X(i,j)=\overline{X(j,i)}=\cases Z_{ij},& i<j,\\ Y_i,&i=j.\endcases
$$

In the case when all the random variables $Y_i,Z_{ij}$ (or $Y_i,\Re Z_{ij},\Im Z_{ij}$
in the Hermitian case) are Gaussian, the Wigner matrix is said to belong to the
Gaussian Orthogonal Ensemble (GOE) in the real case, and Gaussian Unitary
Ensemble (GUE) in the Hermitian case.

For any finite set $B\subset \{1,2,\dots\}$ we denote by $X(B)$ the $|B|\times |B|$
submatrix of the (real symmetric or Hermitian) Wigner matrix $X$ formed by the
intersections of the rows and columns of $X$ marked by elements of $B$. Clearly,
the distribution of $X(B)$ depends only on $|B|$.

Traditionally one encodes the real symmetric and the Hermitian cases by a parameter
$\beta$ that takes value 1 for GOE and value 2 for GUE.

\subhead The height function
\endsubhead
Let $A=\{a_n\}_{n\ge 1}$ be an arbitrary sequence of pairwise distinct natural
numbers. The {\it height function\/} $H_A$ associated to $A$ and a Wigner matrix
$X$ is a random integer-valued function on $\R\times \R_{\ge 1}$ defined by
$$
H_A(x,y)=\sqrt{\frac{\beta\pi}{2}}\,\bigl\{\text{the number of eigenvalues of
}X(\{a_1,\dots,a_{[y]}\})\text{ that
are }\ge x\bigr\}.
$$
The convenience of the constant prefactor $\sqrt{{\beta\pi}/{2}}$ will be
evident shortly.

\subhead Good families of sequences
\endsubhead
In what follows $L>0$ is a large parameter.

Let $\{A_i\}_{i\in I}$ be a family of sequences of pairwise distinct natural numbers.
Assume they all depend on $L$. Denote
$$
A_i=\{a_{i,n}\}_{n\ge 1}, \qquad A_{i,m}=\{a_{i,1},\dots,a_{i,m}\},\qquad i\in I,\quad
m\in\N.
$$

We say that  $\{A_i\}_{i\in I}$ is a {\it good family\/} if for any
$i,j\in I$ and $x,y\in\R_{>0}$
there exists a limit
$$
\alpha(i,x;j,y)=\lim_{L\to\infty}\frac{\left|A_{i,[xL]}\cap A_{j,[yL]}\right|}
{L}\,.
$$
Here is an example of a good family: $I=\{1,2,3,4\}$ and
$$
a_{1,n}=n,\qquad a_{2,n}=2n,\qquad a_{3,n}=2n+1,\qquad a_{4,n}=\cases n+L, & n\le L,\\
n-L,& L<n\le 2L,\\ n,&n>2L.\endcases
$$
Note, however, that the index set $I$ does not have to be finite.

\subhead Correlated Gaussian Free Fields
\endsubhead
Let $\{A_i\}_{i\in I}$ be a good family of sequences as above. Take a family of
copies of the
upper half-plane $\H=\{z\in\C\mid \Im z>0\}$ indexed by $I$ and consider their union
$$
\H(I)=\bigcup_{i\in I} \H_i.
$$

Introduce a function $C:\H(I)\times\H(I)\to \R\cup\{-\infty\}$ via
$$
C_{ij}(z,w)=\frac 1{2\pi}\ln\left|\frac{\alpha(i,|z|^2;j,|w|^2)-zw}
{\alpha(i,|z|^2;j,|w|^2)-z\overline{w}}\right|,\qquad i,j\in I,\quad
z\in\H_i,\ w\in\H_j,
$$
where $\alpha(\,\cdot\,)$ is as above. Note that for $i=j$
$$
C_{ii}(z,w)=\frac 1{2\pi}\ln\left|\frac{\min(|z|^2,|w|^2)-zw}
{\min(|z|^2,|w|^2)-z\overline{w}}\right|=-\frac 1{2\pi}\ln\left|\frac{z-w}
{z-\overline{w}}\right|
$$
is the Green function for the Laplace operator on $\H$ with Dirichlet boundary
conditions.
\proclaim{Proposition 1} For any good family of sequences as above, there exists a
generalized Gaussian process on $\H(I)$ with the covaraince kernel $C(z,w)$ as above.
More exactly, for any finite family of test functions $f_m(z)\in C_0(\H_{i_m})$ and
$i_1,\dots,i_M\in I$, the covariance matrix
$$
cov(f_k,f_l)=\int_\H\int_\H f_k(z)f_l(w) C_{i_ki_l}(z,w) \,dzd\bar{z}\,dwd\bar{w},\qquad
k,l=1,\dots,M,
$$
is positive-definite.
\endproclaim

Denote the resulting generalized Gaussian process by $\Cal G_{\{A_i\}_{i\in
I}}$.

The proof of Propositon 1 will be given later.

\subhead Complex structure
\endsubhead
Let $A$ be a sequence of pairwise distinct integers. The height function $H_A$
is naturally defined on $\R\times \R_{\ge 1}$. Having the large parameter $L$, we would
like to scale $(x,y)\mapsto (L^{-\frac 12}x,L^{-1}y)$, which lands us in $\R\times
\R_{> 0}$.

Wigner's semicircle law implies that with $L\gg 1$, $x\sim L^\frac 12$, $y\sim L$,
after rescaling with overwhelming probability
the eigenvalues (or, equivalently, the places of growth of the height function in $x$-direction) 
are
concentrated in the domain
$$
\bigl\{(x,y)\in\R\times\R_{>0}\mid -2\sqrt{y}\le x\le 2\sqrt{y}\bigr\}.
$$

Let us identify the interior of this domain with $\H$ via the map
$$
\Omega:(x,y)\mapsto\frac x2+ i\sqrt{y-\left(\frac x2\right)^2}.
$$
Its inverse has the form
$$
\Omega^{-1}(z)=(x(z),y(z))=(2\Re(z),|z|^2).
$$
Note that this map sends the boundary of the domain to the real line.

Thanks to $\Omega$ we can now speak of the height function $H_A$ as being defined on
$\H$; we will use the notation
$$
H_A^\Omega(z)=H_A(L^{\frac 12}x(z), Ly(z)),\qquad z\in\H.
$$
Note that we have incorporated rescaling in this definition.

\subhead{Main result}\endsubhead Let $X$ be a (real symmetric or Hermitian) Wigner matrix.
Let $\{A_i\}_{i\in I}$ be a good family of
sequences. We argue that the collection of the centralized random height functions
$$
H_{A_i}^\Omega(z_i)-\E H_{A_i}^\Omega(z_i),\qquad i\in I,\quad z_i\in\H_i,
$$
viewed as distributions, converges as $L\to\infty$ to the generalized Gaussian process
$\Cal G_{\{A_i\}_{i\in I}}$.

One needs to verify the convergence on a suitable set of test functions.
The exact statement that we prove is the following.

\proclaim{Theorem 2}
Pick $i\in I$, $y>0$, and $k\in \Z_{\ge 0}$. Define a moment of the random
height function by
$$
M_{i,y,k}=\int_{-\infty}^{+\infty} x^k \bigl(H_{A_i}(L^{\frac 12}x,Ly)
-\E H_{A_i}(L^{\frac
12}x,Ly)\bigr)dx.
$$
Then as $L\to\infty$, these moments converge, in the sense of finite dimensional
distributions, to the moments of $\Cal G_{\{A_i\}_{i\in I}}$ defined as
$$
\Cal M_{i,y,k} = \int_{z\in\H_i,|z|^2=y} (x(z))^k\, \Cal G_{\{A_i\}_{i\in I}}(z)\,
\frac{dx(z)}{dz}\,dz.
$$
\endproclaim

\subhead Moments as traces
\endsubhead
Let us rescale the variable $x=L^{-\frac 12}u$ in the definition of $M_{i,y,k}$
and then integrate by parts. Since the derivative of the height function
$H_{A_i}(u,[Ly])$ in $u$ is
$$
\frac d{du} H_{A_i}(u,[Ly])=-\sqrt{\frac{\beta\pi}2} \sum_{s=1}^{[Ly]}
\delta(u-\lambda_s),
$$
where $\{\lambda_s\}_{1\le s\le[Ly]}$ are the eigenvalues of $X(A_{i,[Ly]})$, we obtain
$$
\multline
M_{i,y,k}={L^{-\frac{k+1}2}} \sqrt{\frac{\beta\pi}2} \left(\sum_{s=1}^{[Ly]}
\frac{\lambda_s^{k+1}}{k+1}-\E\sum_{s=1}^{[Ly]}
\frac{\lambda_s^{k+1}}{k+1}\right)\\=
\frac{L^{-\frac{k+1}2}} {k+1}\sqrt{\frac{\beta\pi}2}
\left(\tr\bigl(X(A_{i,[Ly]})^{k+1}\bigr)-
\E \tr\bigl(X(A_{i,[Ly]})^{k+1}\bigr)\right).
\endmultline
$$

We can now reformulate the statement of Theorem 2 as follows.

\proclaim{Theorem 2'} Let $X$ be a Wigner matrix. Let $k_1,\dots,k_m\ge 1$ be integers,
and let $B_1,\dots,B_m$ be subsets of $\N$ dependent on the large parameter $L$
such that there exists limits
$$
b_p=\lim_{L\to\infty} \frac{|B_p|}{L}>0,\qquad
c_{pq}=\lim_{L\to\infty}\frac{|B_p\cap B_q|}{L}\,,\qquad p,q=1,\dots,m.
$$
Then the $m$-dimensional random vector
$$
\left(L^{-\frac {k_p}2}\biggl(\tr\bigl(X(B_p)^{k_p}\bigr)-
\E \tr\bigl(X(B_p)^{k_p}\bigr)\biggr)\right)_{p=1}^m
\tag 1
$$
converges (in distribution and with all moments) to the zero mean $m$-dimensional Gaussian
random variable $(\xi_p)_{p=1}^m$ with the covariance
$$
\multline
\E \xi_p\xi_q\\=\frac{2k_pk_q}{\beta \pi}
\oint\limits_{\Sb |z|^2=b_p\\ \Im z>0\endSb} \oint\limits_{\Sb |w|^2=b_q\\ \Im
w>0\endSb}
(x(z))^{k_p-1}(x(w))^{k_q-1}
\frac 1{2\pi}\ln\left|\frac{c_{pq}-zw}
{c_{pq}-z\overline{w}}\right|\,
\frac{dx(z)}{dz}\frac{dx(w)}{dw}\,dzdw.
\endmultline
\tag 2
$$
\endproclaim

\subhead Proof of Theorem 2'
\endsubhead
The argument closely follows that given in Section 2.1.7 of
\cite{AGZ} in the case of one set $B_j\equiv B$. One proves the
convergence of moments, which is sufficient to also claim the convergence
in distribution for Gaussian limits.

Any joint moment of the coordinates of \thetag{1} is written as a finite combination of
contributions corresponding to suitably defined graphs that are in their turn associated
to words. The only difference of the multi-set case with the one-set case is that one
needs to keep track of the {\it alphabets\/} these words are built from: A word
corresponding to coordinate number $p$ of \thetag{1} would have to be built from the
alphabet that coincides with the set $B_p$. Equivalently, the corresponding graphs
will have their vertices labeled by elements of $B_p$.

Since all sizes $|B_p|$ have order $L$, and $|B_1\cup\dots\cup B_m|=O(L)$, the estimate
showing that all contributions not coming from matchings are negligible (Lemma 2.1.34
in \cite{AGZ}) carries over without difficulty. It only remains to compute the
covariance.

For real symmetric Wigner matrices in the one-set case the limits of the
variances of the coordinates of \thetag{1} are given by (2.1.44) in \cite{AGZ}. It
reads (with $k=k_p$ for a $p$ between 1 and $m$)
$$
2k^2 C^2_{\frac{k-1}2}+k^2 C^2_{\frac k2} +\sum_{r=3}^\infty \frac{2k^2}r\left(
\sum_{\Sb k_i\ge 0\\ 2\sum_{i=1}^r k_i=k-r\endSb}\prod_{i=1}^r C_{k_i}\right)^2,
\tag 3
$$
where $\{C_k\}_{k\ge 1}$ are the Catalan numbers, and we assume $C_a=0$ unless
$a\in\{0,1,2,\dots\}$. The Catalan number $C_k$ counts the number of rooted planar
trees with $k$ edges, and different terms of \thetag{3} have the following
interpretation (see \cite{AGZ} for detailed explanations):

$\bullet$ The first term comes from two trees with $(k-1)/2$ edges each that hang
from a common vertex; the factor $k^2$ originates from choices of certain
starting points on each tree united with the common vertex, and the extra 2 is
actually $\E Y_1^2$.

$\bullet$ The second term comes from two trees with $k/2$ edges each
that are glued along one edge. There are $k/2$ choices of this edge for each
of the trees, there is an additional $2=\E Z_{12}^4-1$, and another addional
2 responsible of the choice of the orientation of the gluing.

$\bullet$ The third term comes from two graphs each of which is a cycle of length
$r$ with pendant trees hanging off each of the vertices of the cycle; the total number
of edges in the extra trees being $(k-r)/2$ (this must be an integer). As for the
first term, there is an extra $k^2=k\cdot k$ coming from the choice of the starting
points and also an extra 2 for the choice of the gluing orientation along the cycle.

For each of the three terms the total number of vertices in the resulting graph is equal to $k$,
and if one labels each vertex with a letter from an alphabet of cardinality $|B|$
this would yield a factor of
$$
|B|(|B|-1)\cdots (|B|-k+1)=|B|^k+O(|B|^{k-1}).
$$
Normalization by $|B|^k$ yields \thetag{3}.

In the general case, in order to evaluate the covariance
$$
L^{-\frac {k_p+k_q}2}\E\left[\left(\tr\bigl(X(B_p)^{k_p}\bigr)-
\E \tr\bigl(X(B_p)^{k_p}\bigr)\right)\left(\tr\bigl(X(B_q)^{k_q}\bigr)-
\E \tr\bigl(X(B_q)^{k_q}\bigr)\right)\right]
\tag 4
$$
in the limit, we need to employ the same graph counting, except for the two graphs
being glued now correspond to different values $k_p$ and $k_q$ of $k$, and their
vertices are marked by letters of different alphabets $B_p$ and $B_q$. 

$\bullet$ The first term gives
$2k_pk_q C_{\frac{k_p-1}2}C_{\frac{k_q-1}2}$ for the graph counting, and an extra
$$
|B_p\cap B_q|\cdot(|B_p|-1)(|B_p|-2)\cdots (|B_p|-\tfrac{k_p+1}2)\cdot(|B_q|-1)(|B_q|-2)\cdots
(|B_q|-\tfrac{k_q+1}2)
$$
for the vertex labeling (the factor $|B_p\cap B_q|$ comes from the only
common vertex). Normalized by $L^{-\frac {k_p+k_q}2}$ this yields
$$
2k_pk_q C_{\frac{k_p-1}2}C_{\frac{k_q-1}2} c_{pq}b_p^{\frac{k_p-1}{2}}b_q^{\frac{k_q-1}{2}}.
$$

$\bullet$ The second term has $k_pk_qC_{\frac {k_p}2} C_{\frac {k_q}2}$ from the graph counting
and $c_{pq}^2b_p^{\frac{k_p}{2}-1}b_q^{\frac{k_q}{2}-1}$ from the label counting; a
total of
$$
k_pk_qC_{\frac {k_p}2} C_{\frac {k_q}2}c_{pq}^2b_p^{\frac{k_p}{2}-1}b_q^{\frac{k_q}{2}-1}
$$

$\bullet$ For the third term in the same way we obtain
$$
\sum_{r=3}^\infty \frac{2k_pk_q}r\left(
\sum_{\Sb s_i\ge 0\\ 2\sum_{i=1}^r s_i=k_p-r\endSb}\prod_{i=1}^r C_{s_i}\right)
\left(\sum_{\Sb t_i\ge 0\\ 2\sum_{i=1}^r t_i=k_q-r\endSb}\prod_{i=1}^r C_{t_i}\right)
c_{pq}^r b_p^{\frac{k_p-r}2} b_q^{\frac{k_q-r}2}
$$
Thus, the asymptotic value of the covariance \thetag{4} is
$$
\multline
2k_pk_q C_{\frac{k_p-1}2}C_{\frac{k_q-1}2}  c_{pq}b_p^{\frac{k_p-1}{2}}b_q^{\frac{k_q-1}{2}}+
k_pk_qC_{\frac {k_p}2} C_{\frac {k_q}2}c_{pq}^2b_p^{\frac{k_p}{2}-1}b_q^{\frac{k_q}{2}-1}
\\+\sum_{r=3}^\infty \frac{2k_pk_q}r\left(
\sum_{\Sb s_i\ge 0\\ 2\sum_{i=1}^r s_i=k_p-r\endSb}\prod_{i=1}^r C_{s_i}\right)
\left(\sum_{\Sb t_i\ge 0\\ 2\sum_{i=1}^r t_i=k_q-r\endSb}\prod_{i=1}^r C_{t_i}\right)
c_{pq}^r b_p^{\frac{k_p-r}2} b_q^{\frac{k_q-r}2}
\endmultline
$$

We now use the fact that for any $S=0,1,2,\dots$
$$
\sum_{\Sb s_i\ge 0\\ \sum_{i=1}^r s_i=S\endSb}\prod_{i=1}^r
C_{s_i}=\binom{2S+r}{S}\frac{r}{2S+r},
$$
see (5.70) in \cite{GKP}. This allows us to rewrite the asymptotic
covariance in terms of binomial coefficients:
$$
\multline
2\binom{k_p}{(k_p-1)/2}\binom{k_q}{(k_q-1)/2}c_{pq}b_p^{\frac{k_p-1}{2}}b_q^{\frac{k_q-1}{2}}
\\+
4\binom{k_p}{k_p/2-1}\binom{k_q}{k_q/2-1}
c_{pq}^2b_p^{\frac{k_p-2}{2}}b_q^{\frac{k_q-2}{2}}
\\+\sum_{r=3}^\infty 2r \binom{k_p}{(k_p-r)/2}\binom{k_q}{(k_q-r)/2}
c_{pq}^r b_p^{\frac{k_p-r}2} b_q^{\frac{k_q-r}2}
\\=\sum_{r=1}^\infty 2r \binom{k_p}{(k_p-r)/2}\binom{k_q}{(k_q-r)/2}
c_{pq}^r b_p^{\frac{k_p-r}2} b_q^{\frac{k_q-r}2}
\endmultline
$$

Using the binomial theorem, we can write this expression as a double contour
integral
$$
\frac 2{(2\pi
i)^2}\iint\limits_{const_1=|z|<|w|=const_2}\left(z+\frac{b_p}z\right)^{k_p}
\left(w+\frac{b_q}w\right)^{k_q}\frac{c_{pq}}{b_p}\frac{dzdw}{(\frac{c_{pq}}{b_p}z-w)^2}\,.
\tag 5
$$

Consider the right-hand side of \thetag{2} and assume that $|z|^2=b_p<b_q=|w|^2$.
Observe that
$$
\multline
2\ln\left|\frac{c_{pq}-zw}
{c_{pq}-z\overline{w}}\right|=-2\ln\left|\frac{\frac{c_{pq}}{b_p}z-w}
{\frac{c_{pq}}{b_p}\bar{z}-\overline{w}}\right|
\\\phantom{aaaaa}=-\ln\left(\frac{c_{pq}}{b_p}z-w\right)
+\ln\left(\frac{c_{pq}}{b_p}z-\bar{w}\right)+\ln\left(\frac{c_{pq}}{b_p}\bar{z}-{w}\right)
-\ln\left(\frac{c_{pq}}{b_p}\bar{z}-\bar{w}\right).
\endmultline
$$
This allows us to rewrite the right-hand side of \thetag{2} as a double contour integral over
complete circles in the form
$$
-\frac{k_pk_q}{2\beta \pi^2}
\oint\limits_{|z|^2=b_p} \oint\limits_{|w|^2=b_q}
(x(z))^{k_p-1}(x(w))^{k_q-1}
\ln\left(\frac{c_{pq}}{b_p}z-w\right)
\frac{dx(z)}{dz}\frac{dx(w)}{dw}\,dzdw.
$$

Recalling that $\beta=1$ and noting that
$$
k_p(x(z))^{k_p-1}\frac{dx(z)}{dz}=\frac{d(x(z))^{k_p}}{dz}\,,\quad
k_q(x(w))^{k_q-1}\frac{dx(w)}{dw}=\frac{d(x(w))^{k_q}}{dw}\,,
$$
we integrate by parts in $z$ and $w$ and recover \thetag{5}. The proof for for
$b_p=b_q$ is obtained by continuity of both sides, and to see that the needed
identity holds for $b_p>b_q$ it suffices to observe that both sides are symmetric in
$p$ and $q$.

The argument in the case of Hermitian Wigner matrices is exactly the same, except
in the combinatorial part for the first term the factor 2 is missing due to the
change in $\E Y_1^2$, in the second term 2 is missing due to the change in $\E
|Z_{12}|^4$, and in the third term 2 is missing because there is no choice in the
orientation of two $r$-cycles that are being glued together. \qed

\subhead Proof of Proposition 1
\endsubhead We need to show that for any complex numbers $\{u_k\}_{k=1}^M$
$$
\sum_{k,l=1}^M u_k\overline{u_l}\int_\H\int_\H f_k(z)f_l(w) C_{i_ki_l}(z,w)
\,dzd\bar{z}\,dwd\bar{w}\ge 0.
$$
We can approximate the integration over
the two-dimensional domains by finite sums of one-dimensional
integrals over semi-circles of the form $|z|=const$. On each semi-circle we further
uniformly approximate the (continuous) integrand by a polynomial
in $\Re(z)$. Finally, for the polynomials the nonnegativity follows from
Theorem 2'.\qed

\subhead Chebyshev polynomials
\endsubhead One way to describe the limiting covariance structure in the one-matrix
case is to show that traces of the Chebyshev polynomials of the matrix are asymptotically
independent, see \cite{J}. A similar effect takes place for submatrices as well.

For $n=0,1,2,\dots$ let $T_n(x)$ be the $n$th degree Chebyshev polynomial of the first kind:
$$
T_n(x)=\cos(n\arccos x),\qquad T_n(\cos(x))=\cos(nx).
$$
For any $a>0$, let $T_n^a(x)=T_n(\frac xa)$ be the rescaled version of $T_n$.

\proclaim{Proposition 3} In the assumptions of Theorem 2', for any $p,q=1,\dots,m$
$$
\multline
\lim_{L\to\infty}\E\Biggl[\left(\tr\bigl(T_{k_p}^{2\sqrt{b_pL^{k_p}}}(X(B_p))\bigr)-
\E \tr\bigl(T_{k_p}^{2\sqrt{b_pL^{k_p}}}(X(B_p))\bigr)\right)\\ \times \left(\tr\bigl(T_{k_q}^{2\sqrt{b_qL^{k_q}}}(X(B_q))\bigr)-
\E \tr\bigl(T_{k_p}^{2\sqrt{b_pL^{k_q}}}(X(B_q))\bigr)\right)\Biggr]
=\delta_{k_pk_q}\,\frac{k_p}{2\beta} \left(\frac{c_{pq}}{\sqrt{b_pb_q}}\right)^{k_p}.
\endmultline
$$
\endproclaim
\demo{Proof}
Using \thetag{5}  and assuming $b_p<b_q$ we obtain
$$
\multline
\lim_{L\to\infty}\E\Biggl[\left(\tr\bigl(T_{k_p}^{2\sqrt{b_pL^{k_p}}}(X(B_p))\bigr)-
\E \tr\bigl(T_{k_p}^{2\sqrt{b_pL^{k_p}}}(X(B_p))\bigr)\right)\\ \times \left(\tr\bigl(T_{k_q}^{2\sqrt{b_qL^{k_q}}}(X(B_q))\bigr)-
\E \tr\bigl(T_{k_p}^{2\sqrt{b_pL^{k_q}}}(X(B_q))\bigr)\right)\Biggr]\\
=
\dfrac 2{\beta(2\pi
i)^2}\displaystyle\iint\limits_{b_p=|z|<|w|=b_q}T_{k_p}(\cos(\arg(z))
T_{k_q}(\cos(\arg(w))\dfrac{c_{pq}}{b_p}\dfrac{dzdw}{(\frac{c_{pq}}{b_p}z-w)^2}\\
=\dfrac 1{2\beta(2\pi
i)^2}\displaystyle\iint\limits_{b_p=|z|<|w|=b_q}\Biggl(\biggl(\frac{z}{\sqrt{b_p}}\biggr)^{k_p}+
\biggl(\frac{\sqrt{b_p}}{z}\biggr)^{k_p}\Biggr) \Biggl(\biggl(\frac{w}{\sqrt{b_q}}\biggr)^{k_q}+
\biggl(\frac{\sqrt{b_q}}{w}\biggr)^{k_q}\Biggr)\\ \times
\dfrac{c_{pq}}{b_p}\dfrac{dzdw}{(\frac{c_{pq}}{b_p}z-w)^2}.
\endmultline
$$
Writing $(\frac{c_{pq}}{b_p}z-w)^{-2}$ as a series in $z/w$ we arrive at the result. 
Continuity and symmetry of both sides of the limiting relation removes the assumption $b_p<b_q$.\qed

\end